\documentclass[a4paper,12pt]{article}
\usepackage{times, url}
\usepackage{graphicx,amsmath,amssymb,amsthm,amsfonts,bm,float,cite}
\newtheorem{theorem}{Theorem}[section]

\newtheorem{definition}{Definition}[section]
\newtheorem{proposition}{Proposition}[section]

\usepackage[none]{hyphenat}

\usepackage[center]{caption}
\renewcommand{\thefootnote}{\fnsymbol{footnote}}

\newcommand\blfootnote[1]{%
  \begingroup
  \renewcommand\thefootnote{}\footnote{#1}%
  \addtocounter{footnote}{-1}%
  \endgroup
}

\usepackage[bottom,symbol]{footmisc}
\DeclareFontEncoding{TS1}{}{}
\DeclareFontSubstitution{TS1}{cmr}{m}{n}
\DeclareTextSymbol{\tcrp}{TS1}{'251}
\DeclareTextSymbolDefault{\tcrp}{TS1}
\newcommand{\licen}{ \phantom{.} \\  \hspace{-10mm} {\small  \fbox{
	\begin{tabular}{p{0.105\textwidth} p{0.82\textwidth}}
	  \raisebox{-22pt}{\includegraphics[scale=0.34]{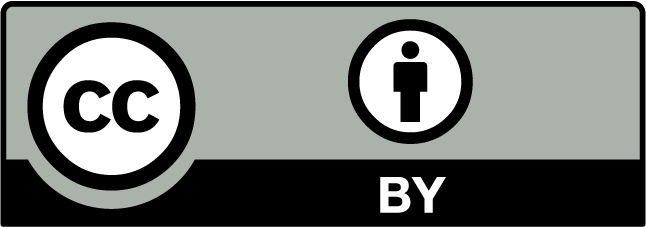}}	& {
			Copyright \tcrp ~2023 by the Authors. This is an Open Access paper distributed under the terms and conditions of the Creative Commons Attribution 4.0 International License (CC BY 4.0).} { {\footnotesize \ttfamily https://creativecommons.org/licenses/by/4.0/}}    \vspace{-0.5mm}
	\end{tabular} }
} \vspace{-10mm}
}

\begin{document}
\setcounter{page}{1}

\thispagestyle{empty} 

\noindent \textbf{Notes on Number Theory and Discrete
Mathematics \newline
Print ISSN 1310--5132, Online ISSN 2367--8275 \newline
XXXX, Volume XX, Number X, XXX--XXX \newline
DOI: 10.7546/nntdm.XXXX}  
\vspace{11mm}

\begin{center}
{\LARGE \bf  An alternative proof of Sylvester's theorem \\[4mm] and variations for more primes} 
\vspace{8mm}

{\Large \bf Steven Brown$^1$}
\vspace{3mm}

$^1$ e-mail: \url{steven.brown.math@gmail.com}
\vspace{2mm}

\end{center}
\vspace{9mm}  \blfootnote{\licen}

\noindent
\noindent {\bf Received:}  13 February 2023 \hfill {\bf Revised:} DD Month XXXX \\
{\bf Accepted:} DD Month XXXX  \hfill {\bf Online First:} DD Month XXXX \\[4mm] 
{\bf Abstract:} This document presents an alternative proof of Sylvester's theorem stating that "the product of $n$ consecutive numbers strictly greater than $n$ is divisible by a prime strictly greater than $n$". In addition, the paper through theorem \ref{thm:gen_ths} proposes stronger versions of Sylvester's theorem and Bertrand's postulate with more primes, as well as an approach to get more results in this direction. 
\vspace{5mm}

\section{Introduction}

A theorem from Sylvester states that the product of $n$ consecutive numbers strictly greater than $n$ is divisible by a prime strictly greater than $n$. \cite{sylvester1892arithmetical}. The purpose of this paper is to demonstrate a theorem (theorem \ref{thm:gen_ths}) which allows
\begin{itemize}
\item to prove Sylvester's theorem
\item to prove stronger versions of the theorem; for example that for $n$ greater than 6, the product of $n$ consecutive numbers strictly greater than $n$ is divisible by two primes strictly greater than $n$
\item to prove a stronger version of Bertrand's postulate
\end{itemize}

The next section "Preliminary works" gives simple results that can be read as and when necessary when referred to from the following sections. The section "Generalization of Sylvester's theorem" gives a theorem and its proof as a basis for stronger results with more than one prime as in Sylvester's theorem. The section "Some application of theorem \ref{thm:gen_ths}" details some of the results mentioned before.
In the rest of the document most variables are integers unless they are seen as real values when considering functions and derivatives. To finalize my proofs I developed python code available on demand. In the rest of the document $\log$ means natural logarithm.

\section{Preliminary works}

\begin{proposition}\label{prop:g_decreases} Let $m$ and $n$ be two integers such that $m\geq n\geq 2$. Let $g$ be the function:
\begin{equation*}
g(x):=\frac{\log(m+x)}{\log(n+x)}
\end{equation*}
$g$ is decreasing for $x\geq 0$.
\end{proposition}

\begin{proof}
The derivative of $g$ with respect to $x$ is
\begin{align*}
g'(x)&=\frac{\frac{\log(n+x)}{m+x}-\frac{\log(m+x)}{n+x}}{{\log(n+x)}^2} \\
&= \frac{(n+x)\log(n+x)-(m+x)\log(m+x)}{(m+x)(n+x){\log(n+x)}^2}
\end{align*}
The derivative is negative for $x\geq 0$, therefore $g$ is a decreasing function of $x$ for $x\geq 0$.
\end{proof}

\begin{definition}\label{def:Pmn}
Let $m$ and $n$ be two integers such that $m\geq n\geq 2$. Let $P(m,n)$ be the integer
\begin{equation*}
P(m,n):=\prod_{i=1}^n (m+i)
\end{equation*}
the product of $n$ consecutive numbers which first term is $m+1$.
\end{definition}

For an integer $a$ and a prime $q$ we write $v_q(a)$ the greatest exponent $b$ such that $q^b$ divides $a$.

\begin{proposition}\label{prop:boundsLegendre} With $P(m,n)$ as defined in definition \ref{def:Pmn} and for any prime $q$
\begin{equation*} 
v_q(n!)\leq v_q(P(m,n)) \leq v_q(n!) + \left\lfloor \frac{\log(m+n)}{\log(q)} \right\rfloor
\end{equation*}
\end{proposition}

\begin{proof}
Referring to \cite{legendre1830theorie},
for any prime $q$, the exponent of the largest power of $q$ that divides $P(m,n)$ is given by

\begin{equation*}
v_q\left(P(m,n)\right)=v_q\left((m+n)!\right)-v_q\left(n!\right)
= \sum_{i=1}^\infty \left( \left\lfloor \frac{m+n}{q^i}\right\rfloor - \left\lfloor \frac{m}{q^i}\right\rfloor \right)
\end{equation*}

From general properties of the floor function, for any real numbers $x$ and $y$ we have
\begin{equation*}
\left\lfloor x \right\rfloor + \left\lfloor y \right\rfloor \leq \left\lfloor x + y \right\rfloor \leq \left\lfloor x \right\rfloor + \left\lfloor y \right\rfloor + 1
\end{equation*}

Therefore for any exponent $i\geq 1$

\begin{equation}\label{Legendre_bounds}
\left\lfloor \frac{n}{q^i}\right\rfloor \leq \left\lfloor \frac{m+n}{q^i}\right\rfloor - \left\lfloor \frac{m}{q^i}\right\rfloor \leq \left\lfloor \frac{n}{q^i}\right\rfloor + 1
\end{equation}

When $i$ is such that $\left\lfloor \frac{m+n}{q^i}\right\rfloor = 0$, we must have $\left\lfloor \frac{m}{q^i}\right\rfloor = 0$ and $\left\lfloor \frac{n}{q^i}\right\rfloor = 0$. \\
Now, the infinite sum of all the inequalities \ref{Legendre_bounds} is equal to the finite sum limited to indexes $i$ not greater than $\frac{\log(m+n)}{\log(q)}$. When then have:

\begin{equation*}
v_q(n!)\leq v_q(P(m,n)) \leq v_q(n!) + \left\lfloor \frac{\log(m+n)}{\log(q)} \right\rfloor
\end{equation*}

\end{proof}

\begin{proposition}\label{prop:bounds_cbc}
If $n\geq 2$ the following bounds hold true for the central binomial coefficient
\begin{equation*}
\forall n\geq 1 \quad \frac{4^n}{\sqrt{(\pi (n+\frac{1}{2}))}}\leq \binom{2n}{n}\leq \frac{4^n}{\sqrt{(\pi n)}}
\end{equation*}
\end{proposition}

\begin{proof} This is a classical result derived from Wallis integral. Let $n$ be a positive integer.
\begin{equation*}
W_n = \int_0^{\frac{\pi}{2}}\left(\cos(x)\right)^n dx
\end{equation*}
$W_0=\frac{\pi}{2}$ and $W_1=1$. For $n\geq 2$ integrating by parts leads to the following relation $nW_n=(n-1)W_{n-2}$. We get to the the formulas $W_{2n}=\frac{\binom{2n}{n}}{4^n}\frac{\pi}{2}$ and $W_{2n+1}=\frac{4^n}{(2n+1)\binom{2n}{n}}$. Now we observe that for $x\in]0;\frac{\pi}{2}[$, $\cos(x)\in ]0;1[$ and we get to
\begin{equation*}
W_{2n+1}<W_{2n}<W_{2n-1}=\frac{2n+1}{2n}W_{2n+1}
\end{equation*}
Which leads to
\begin{equation*}
\frac{4^n}{(2n+1)\binom{2n}{n}}<\frac{\binom{2n}{n}}{4^n}\frac{\pi}{2}<\frac{4^n}{2n\binom{2n}{n}}
\end{equation*}

And the result follows.
\end{proof}

\begin{definition}\label{def:kmn} For $m\geq n \geq 2$, Let $k(m,n)$ be the following function
\begin{equation*}
k(m,n):= \frac{\sum_{i=1}^n\log\left(1+\frac{m}{i}\right)}{\log(m+n)}
\end{equation*}
\end{definition}

\begin{proposition}\label{prop:kan_sg_kbn} Let $a$, $b$ and $n$ be integers such that $a > b\geq n \geq 2$ then $k(a,n)>k(b,n)$
\end{proposition}

\begin{proof} Assume that $m>n\geq 2$,
\begin{multline*}
k(a,n)-k(b,n) = \frac{\sum_{i=1}^n\log\left(1+\frac{a}{i}\right)}{\log(a+n)}-\frac{\sum_{i=1}^n\log\left(1+\frac{b}{i}\right)}{\log(b+n)} \\
= \sum_{i=1}^n\frac{\log(b+n)\left(\log(a+i)-\log(i)\right)-\log(a+n)\left(\log(b+i)-\log(i)\right)}{\log(a+n)\log(b+n)} \\
= \sum_{i=1}^n \frac{\log(b+n)\log(a+i)-\log(a+n)\log(b+i)}{\log(a+n)\log(b+n)} + \\
 \sum_{i=1}^n \log(i)\frac{\log(a+n)-\log(b+n)}{\log(a+n)\log(b+n)}
\end{multline*}
The second sum is strictly positive given that $\log(a+n)>\log(b+n)$ and $n\geq 2$ ensures there is at least one term not zero since $\log(2)>0$.
The first sum is also positive. From proposition \ref{prop:g_decreases}, for any $i\leq n$ we must have $g(i)\geq g(n)$. That is
\begin{equation*}
\frac{\log(a+i)}{\log(b+i)} \geq \frac{\log(a+n)}{\log(b+n)}
\end{equation*}
Hence
\begin{equation*}
\log(b+n)\log(a+i)-\log(a+n)\log(b+i)\geq 0
\end{equation*}
Which ensures that the first sum is positive. We then have $k(a,n)>k(b,n)$.
\end{proof}

\begin{definition}\label{def:En}
Let $C$ be the constant $1.25506$ and $E(n)$ the function defined, for $n>1$ by
\begin{equation*}
E(n):=\frac{n}{\log(n)}\left(\log(4)-C-\frac{\log(2)\log(4)}{\log(2n)}\right)
\end{equation*}
\end{definition}

\begin{proposition}\label{prop:En} With the function $E$ as defined in definition \ref{def:En}
\begin{itemize}
\item  if $m>n\geq 3$ then $E(m)>E(n)$
\item $E(n) \underset{+\infty}{\sim} \frac{n}{\log(n)}(\log(4)-C)$
\end{itemize}
\end{proposition}

\begin{proof} Regarding the first point, if $\alpha(n)=\frac{n}{\log(n)}$ and $\beta(n)=\log(4)-C-\frac{\log(2)\log(4)}{\log(2n)}$ are two function defined for $n>1$ we have $E(n)=\alpha(n) \beta(n)$. Both $\alpha$ and $\beta$ are increasing functions of $n$ for $n>1$ as on can see looking at their first derivatives.
\begin{equation*}
\alpha'(n)=\frac{\log(n)-1}{{\log(n)}^2}
\end{equation*}
$\alpha'$ is strictly positive for $n> e$.
\begin{equation*}
\beta'(n)=\frac{\log(2)\log(4)}{2n}
\end{equation*}
$\beta'$ is strictly positive for $n>1$.
Therefore for $n\geq 3$, both $\alpha$ and $\beta$ are strictly increasing functions of $n$. Then for $m>n\geq 3$ we have
\begin{equation*}
E(m)=\alpha(m)\beta(m)>\alpha(n)\beta(n)=E(n)
\end{equation*}

The second point is clear when writing
\begin{equation*}
E(n)=\frac{n}{\log(n)}\left(\log(4)-C\right)\left(1-\frac{1}{\log(2n)}\frac{\log(2)\log(4)}{\log(4)-C}\right)
\end{equation*}
with
\begin{equation*}
\underset{n\rightarrow +\infty}{\lim} 1-\frac{1}{\log(2n)}\frac{\log(2)\log(4)}{\log(4)-C} = 1
\end{equation*}
\end{proof}

\section{Generalization of Sylvester theorem}

\begin{theorem}\label{thm:gen_ths} With the function $E$ as defined in \ref{def:En}, anytime there exists integers $n^\star>1$ and $r\geq 0$ such that $E(n^\star)> r + 1$, then, for any $m \geq n\geq n^\star$, $P(m,n)$ as defined in definition \ref{def:Pmn} is divisible by more than $r$ distinct primes all strictly greater than $n$. That is, the product of $n$ consecutive numbers strictly greater than $n$ is divisible by more than $r$ distinct primes all strictly greater than $n$.
\end{theorem}

\begin{proof}
Let $m$ and $n$ be two integers such that $m\geq n\geq 2$ and $P(m,n)$ be the product of $n$ consecutive numbers which first term is $m+1$ as in definition \ref{def:Pmn}. Let $r\geq 1$ be an integer and let $p_1,\ldots,p_k$ be all the prime numbers lower than $n$, with $k=\pi(n)$, $\pi$ being here the prime counting function. Let's assume as well that there exists 
$r$ distinct primes strictly greater than $n$, say $q_1,\ldots,q_r$, dividing $P(m,n)$. Each prime $q_i$ divides only one term of the sequence because $q_i>n$, therefore its maximum exponent $\beta_i$ in $P(m,n)$ must satisfy
\begin{equation*}
{q_i}^{\beta_i}\leq m+n
\end{equation*}
Now assume that the support of $P(m,n)$ is exactly $\prod_{i=1}^k p_i \prod_{j=1}^r q_j$, then by virtue of the above remark and of proposition \ref{prop:boundsLegendre} we must have

\begin{equation*}
P(m,n) = \prod_{i=1}^k {p_i}^{\alpha_i} \prod_{j=1}^r {q_j}^{\beta_j} \leq (m+n)^r n! \prod_{i=1}^k {p_i}^{\left\lfloor \frac{\log(m+n)}{\log(p_i)}  \right\rfloor} \leq (m+n)^r n! \prod_{i=1}^k {p_i}^{\frac{\log(m+n)}{\log(p_i)}}
\end{equation*}

With $k(m,n)$ as in definition \ref{def:kmn}, the above inequality is tantamount to

\begin{equation}\label{ineq:fund}
k(m,n)=\frac{\log\left(\binom{m+n}{n}\right)}{\log{m+n}}\leq r + \pi(n)
\end{equation}

By way of contradiction, assuming all the conditions above, the denial of this inequality \ref{ineq:fund} implies that there are strictly more than $r$ distinct primes stricly greater than $n$ dividing $P(m,n)$.

Let $C$ be the constant $C=1.25506$ and $b(n)$ the function
\begin{equation*}
b(n):= C\frac{n}{\log(n)}
\end{equation*}
From \cite{rosser1962approximate} the following holds true for $n> 1$
\begin{equation}\label{inequation:rosser}
\pi(n)<b(n)
\end{equation}

Let's find a lower bound to $k(n,n)-b(n)$. From proposition \ref{prop:bounds_cbc} we have

\begin{align*}
k(n,n)-b(n) &\geq \frac{n\log(4)}{\log(2n)} - \frac{1}{2}\frac{\log(\pi(n+\frac{1}{2}))}{\log(2n)} - C\frac{n}{\log(n)} \\ 
&\geq \frac{n\log(4)}{\log(2n)} - C\frac{n}{\log(n)} - 1
\end{align*}

Given that
\begin{equation*}
\frac{\log(4)}{\log(2n)}=\frac{\log(4)}{\log(n)}-\frac{\log(2)\log(4)}{\log(n)\log(2n)}
\end{equation*}
We get to
\begin{equation}\label{ineq:knn}
k(n,n) > \pi(n) + E(n) - 1
\end{equation}

Let $n^\star$ be an integer such that $E(n^\star)> r+1$. Such an integer exists referring to proposition \ref{prop:En}. Now let $m$ and $n$ be integers such that $m\geq n \geq n^\star$, From proposition \ref{prop:En} we have
\begin{equation*}
E(n) \geq E(n^\star)> r+1
\end{equation*}
From inequality \ref{ineq:knn} we must have
\begin{equation*}
k(n,n) > \pi(n) + E(n) - 1 > \pi(n) + r
\end{equation*}
And from proposition \ref{prop:kan_sg_kbn} (with $a=m$ and $b=n$) for $m\geq n$, we must have $k(m,n)\geq k(n,n)$. This leads to
\begin{equation}\label{ineq:denial}
k(m,n) > \pi(n) + r
\end{equation}
This inequality \ref{ineq:denial} is a denial of inequality \ref{ineq:fund}, therefore for any integers $m$ and $n$ satisfying $m\geq n \geq n^\star$ there are strictly more than $r$ distinct primes strictly greater than $n$ dividing $P(m,n)$.
\end{proof}

Note that we could have used a better upper bound of $\frac{1}{2}\frac{\log(\pi(n+\frac{1}{2}))}{\log(2n)}$ in the previous theorem, however for the purpose of the applications of the theorem to follow this seemed to have little added value.

\section{Some implications of theorem \ref{thm:gen_ths} }

\subsection{A proof of Sylvester's theorem}
We have $E(1100)>1$; and the condition of theorem \ref{thm:gen_ths} with $n^\star=1100$ and $r=0$. More than zero primes being at least one prime, by application of the theorem, Sylvester's theorem is satisfied for any $m\geq n \geq 1100$.
We finish the proof with an extensive analysis of the cases when $2\leq n<1100$. With a program, for each $n$ we find $m^\star\geq n$ such that $k(m^\star,n)>\pi(n)$ and, for that particular $n$, the theorem is valid for any $m\geq m^\star$. This means we are left with only a finite number of cases to check and an exhaustive analysis of these cases, e.g. the existence of a prime strictly greater than $n$ in the factorization of $P(m,n)$ ends the proof. The Python program doing these calculations is available on demand.

\subsection{More theorems of that kind}

From proposition \ref{prop:En}, $E(n)$ is equivalent to $\frac{n}{\log(n)}(\log(4)-C)$ as $n$ grows to infinity. That means for any $r\geq 0$ we will find $n^\star$ such that $E(n^\star)> r+1$ and the theorem \ref{thm:gen_ths} will apply. Therefore, for any integers $m\geq n \geq n^\star$ the sequence $m+1,m+2,\ldots ,m+n$ will contain strictly more than $r$ distinct primes strictly greater that $n$. With the help of programming, if the calculation is doable, it is also possible to prove the theorem for lower values of $m$ and $n$ allowing reducing $n^\star$. For example:

\begin{theorem} With $n\geq 6$, the product of $n$ consecutive numbers strictly greater than $n$ is divisible by at least two distinct primes strictly greater than $n$.
\end{theorem}

\begin{proof} Applying theorem \ref{thm:gen_ths} with $E(1411)>2$, that is $n^\star=1411$ and $r=1$ and checking all cases for $6\leq n<1411$. For each $6\leq n<1411$ we find $m^\star \geq n$ such that $k(m^\star ,n)>\pi(n) + 1$ leaving us with only a finite number of cases to check because $l\geq m \geq n$ implies $k(l,n)\geq k(m,n)$ from proposition \ref{prop:kan_sg_kbn}
\end{proof}

\begin{theorem} With $n\geq 9$, the product of $n$ consecutive numbers strictly greater than $n$ is divisible by at least three distinct primes strictly greater than $n$.
\end{theorem}

\begin{proof} The same method with $E(1705)>3$ and being able to check a finite number of cases
\end{proof}

More theorems of than kind can easily be derived from this process.

\subsection{Number of primes in $]n,2n[$}

From Bertrand's postulate, for any integer $n>3$ there always exists at least one prime number $p$ satisfying $n<p<2n-2$. We shall see that this result can be improved with theorem \ref{thm:gen_ths}. 

\begin{theorem} For $n\geq 1100$, the number of primes in $]n,2n[$ is greater than
\begin{equation*}
\left\lfloor \frac{n}{\log(n)}\left(\log(4)-C-\frac{\log(2)\log(4)}{\log(2n)}\right) \right\rfloor
\end{equation*}
\end{theorem}

\begin{proof} For $n\geq 1100$, $E(n)>1$ and $E(n)\geq \left\lfloor E(n) \right\rfloor \geq 1$. Then
\begin{equation*}
E(n) \geq \left( \left\lfloor E(n) \right\rfloor - 1 \right) + 1
\end{equation*}
We are in the conditions of theorem \ref{thm:gen_ths} with $r=\left\lfloor E(n) \right\rfloor - 1\geq 0$. Taking $m=n$ we conclude there are strictly more than $\left\lfloor E(n) \right\rfloor - 1$ primes in $]n;2n[$. which is equivalent to say there are at least $\left\lfloor E(n) \right\rfloor$ primes in $]n;2n[$.
\end{proof}

\begin{theorem} For $n\geq 1123$, the number of primes in $]n,2n[$ is greater than
\begin{equation*}
\left\lfloor \frac{1}{C} \left( \log(4) - C - \frac{\log(2)\log(4)}{\log(2n)}\right)\pi(n) \right\rfloor
\end{equation*}
\end{theorem}

\begin{proof} From \cite{rosser1962approximate} we can write that
\begin{equation*}
E(n)>\frac{\pi(n)}{C} \left( \log(4) - C - \frac{\log(2)\log(4)}{\log(2n)} \right) := \rho(n)
\end{equation*}
Because of the condition $n\geq 1123$, we have $\rho(n) > 1$ and $E(n)>\left\lfloor \rho(n) \right\rfloor \geq 1$.
Then we have
\begin{equation*}
E(n)>\left( \left\lfloor \rho(n) \right\rfloor - 1 \right) + 1
\end{equation*}
We have the conditions to apply theorem \ref{thm:gen_ths} with $r=\left\lfloor \rho(n) \right\rfloor - 1\geq 0$. We take $m=n$ and the conclusion follows.
\end{proof}

For example, with $n^\star = 620\;634$ we have 
\begin{equation*}
\frac{1}{C} \left( \log(4) - C - \frac{\log(2)\log(4)}{\log(2n^\star)}\right)>0.05
\end{equation*}
Given that $n\geq n^\star$ implies $E(n)\geq E(n^\star)$ from proposition \ref{prop:En} we can see that the number of primes in $]n;2n[$ for $n \geq n^\star$ is greater than $\left\lfloor 0.05 \cdot \pi(n) \right\rfloor$.

The maximum ratio that we can find with this method is not greater than:
\begin{equation*}
\frac{\log(4)}{C} - 1 < 0.104565
\end{equation*}
But for any ratio $\theta < \frac{\log(4)}{C} - 1$ we can find $n^\star$ such that for any $n\geq n^\star$, the number of primes in $]n;2n[$ is greater than $\left\lfloor  \theta \pi(n) \right\rfloor$. $\theta$ can be chosen as high as $10\%$.
 
\makeatletter
\renewcommand{\@biblabel}[1]{[#1]\hfill}
\makeatother

\bibliographystyle{plain}

\nocite{*}

\bibliography{citation}

%
%
%
%
%
%
%
\end{document}